\documentclass[preprint,12pt]{elsarticle}
\journal{Discrete Mathematics}

\usepackage[utf8]{inputenc}
\usepackage{amsmath}
\usepackage{amsfonts}
\usepackage{amssymb}
\usepackage{graphicx}
\usepackage[left=2cm,right=2cm,top=2cm,bottom=2cm]{geometry}

\newtheorem{thm}{Theorem}
\newtheorem{lemma}[thm]{Lemma}
\newtheorem{corollary}[thm]{Corollary}

\begin{document}

%% We can define new commands and use them in the paper!
\newcommand{\TT}{T-tetromino}
\newcommand{\TTs}{T-tetrominos}

\begin{frontmatter}
\title{T-Tetrominos in Arithmetic Progression}

\author[emily]{Emily Feller}

\affiliation[emily]{
    organization={Acacia International School},
    addressline={Plot 1222 Zzimwe Road},
    city={Kampala},
    country={Uganda}}
  
\author[rob]{Robert Hochberg}

\affiliation[rob]{
  organization={University of Dallas},
  addressline={1845 E. Northgate Drive}, 
  city={Irving}, 
  state={TX},
  postcode={75062},
  country={USA}}

\begin{abstract}
A famous result of D. Walkup is that an $m\times n$ rectangle may be tiled by \TTs\ if and only if both $m$ and $n$ are multiples of 4. The ``if'' portion may be proved by tiling a $4\times 4$ block, and then copying that block to fill the rectangle; but this leads to regular, periodic tilings. In this paper we investigate how much ``order'' must be present in every tiling of a rectangle by \TTs, where we measure order by length of arithmetic progressions of tiles.
\end{abstract}


\end{frontmatter}

%% Also, style note: We can't put "Feller unit" in the paper, since there is an old, established, universal and inviolable tradition of not naming things after yourself, even if it was your advisor who named it. So in the writeup here, I just called them units. 

\section{Introduction}
There is a decent-sized body of literature characterizing the rectangles that may be tiled by various polyominos, and the \TT\ has received much attention. A well-known result of D. Walkup \cite{Walkup} is that the \TT\ can tile a rectangle if and only if both sides of that rectangle are multiples of 4. This result has been extended in several ways. Mike Reid \cite{ReidKlarner} generalized the \TT\ to a class of $(8n-4)$-ominos, each having Walkup-like conditions on the rectangles that they tile. Authors \cite{KornPak, Merino} have enumerated tilings of $4m\times 4n$ rectangles with $T$-tetrominos, and others \cite{Goddard, Hochberg, Zhan} have  characterized (partial) tilings of rectangles that do not meet Walkup's condition.

%The textures are still difficult to differentiate in black and white.  Is it possible to make different textures for each color? Emily, I think I have addressed this now, but I don't think it is very pretty...

\begin{figure}[b]
\begin{center}
\includegraphics[scale=0.9]{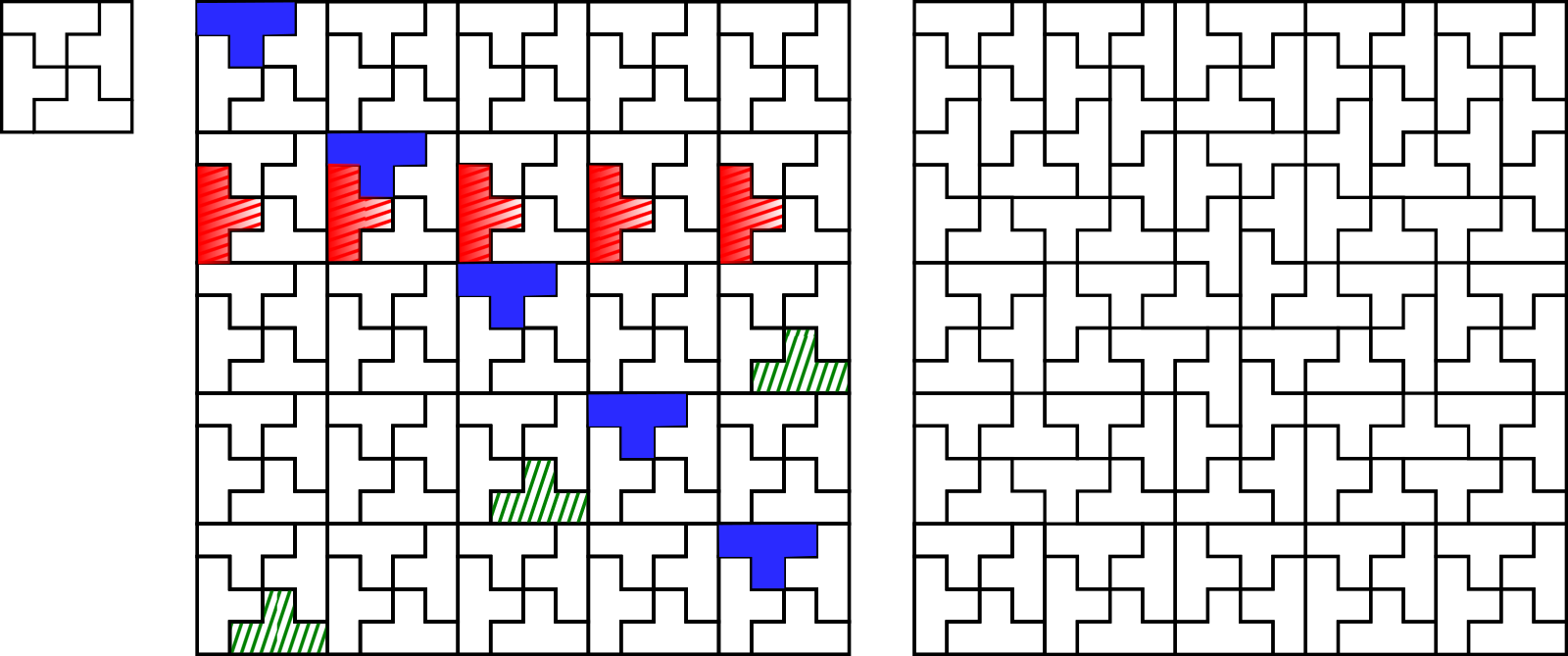}
\end{center}
\caption{(Left) A 4x4 unit. (Middle) A periodic tiling of a $20\times 20$ square. (Right) A tiling without any 3-term AP of tiles.}
\label{fig-4x4}
\end{figure}

Any $4m\times 4n$ rectangle may be tiled by repeating the four-tile unit on the left in Figure \ref{fig-4x4}, yielding an orderly tiling of, for example, a $20\times 20$ rectangle as shown in the middle of that figure. The $20\times 20$ tiling shown on the right in that figure is visually less orderly. We can make this precise by introducing the notion of an {\em arithmetic progression} (AP) of tiles; that is, a sequence $T_1, T_2, \ldots, T_l$ of tiles that are all in the same orientation, and such that the translation from $T_i$ to $T_{i+1}$ is the same for all $i, 1\leq i \leq l-1$. The number $l$ of tiles is called the {\em length} of the AP. The tiling in the middle of Figure \ref{fig-4x4} has several APs of length 5 (two are highlighted, along with one of length three), but the tiling on the right has no AP of length greater than two.

The natural Ramsey-type question to ask is, given a length $l$, must there exist numbers $M$ and $N$ so that whenever $m\geq M$ and $n\geq N$, any tiling of a $4m\times 4n$ rectangle with \TTs\ must contain an AP of tiles of length at least $l$? And from a fixed rectangle's point of view we can ask, given $m$ and $n$, what is the largest length AP of tiles that every $4m\times 4n$ tiling must contain?

In this paper, ``tiling'' will always refer to tiling by non-overlapping \TTs, and tilings of rectangles will always mean completely covering the interior of a rectangle satisfying Walkup's condition. We write $(m, n)\rightarrow l$ to mean that every tiling of an $m\times n$ rectangle must contain an arithmetic progression of tiles of length at least $l$. For a fixed width $w$, a multiple of 4, and $l$, any integer at least 1, we define $T = T_w(l)$ to be the least value of $T$ such that $(w, T)\rightarrow l$.

The rest of the paper is organized as follows: 
In Section \ref{sec-basic} we prove some basic results and make a connection with classical van der Waerden numbers. In Sections \ref{sec-w4} and \ref{sec-chain} we show that for $w$ in $\{4, 8, 12, 16\}$ and $l\geq 3$, $T_w(l)$ is exactly given by the two-color van der Waerden numbers, and that $T_{20}(l)$ is not. We end with a section on exact values and further explorations.
%MMM that good above? "is not"

\section{Basic Results}\label{sec-basic}

\begin{figure}[h]
\begin{center}
\includegraphics[scale=0.6]{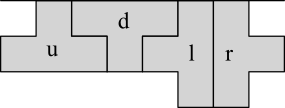}
\end{center}
\caption{Boundary tiling with no arithmetic progression of length 2.}
\label{fig-boundary}
\end{figure}

We begin by considering {\em boundary tilings}, that is, coverings of unit squares along one side of a line, by tiles lying completely on that side of the line, as shown in Figure \ref{fig-boundary}. The letters in the left figure stand for {\em up, down, left} and {\em right}, the orientations of the tiles --- wording and notation we will frequently use.

We start with a simple theorem:

\begin{thm}
Any boundary tiling covering seven consecutive boundary squares must contain an arithmetic progression of length two.
\end{thm}

\noindent{\bf Proof:} If there are two tiles in the same orientation then we are done. But four tiles, one in each orientation, can cover at most six squares along that edge.\hfill$\Box$

\begin{thm}\label{thm-boundary}
For any length $l\geq1$ there exists some boundary length $B$ such that any covering of $B$ consecutive boundary squares must contain an arithmetic progression of length $l$.
\end{thm}

The proof makes use of van der Waerden's theorem. See \cite{GRS}.

\begin{thm}[van der Waerden's Theorem]\label{thm-vdW}
For any integers $r\geq 1, l\geq 1$ there exists an integer $W=W(r,l)$ such that any $r$-coloring of the initial segment of positive integers 
$\{1,2, 3, \ldots, W\}$ contains a monochromatic arithmetic progression of length $l$.
\end{thm}

\noindent{\bf Proof of Theorem \ref{thm-boundary}:} 
Color the boundary squares using the colors $u$, $l$ and $r$ for squares covered by tiles in those three orientations, and $d$1, $d$2 and $d$3 for squares covered by the left, middle or right squares of a $d$ tile. This yields a six-coloring of the boundary squares. If we let $B=W(6, l)$ then there must be a monochromatic AP of length $l$, which corresponds to an AP of Ts of length $l$.
\hfill$\Box$

We will improve this bound following the discussion of the HV Lemma below. As a straightforward corollary we get a Ramsey-type theorem about APs of \TTs\ in tilings of rectangles.

\begin{corollary}
For any length $l\geq 1$ there exists some number $B$ such that any tiling of a rectangle having some side of length at least $B$ must contain an AP of \TTs\ of length $l$.
\end{corollary}

\section{Width 4}\label{sec-w4}

For rectangles of width 4, 8, 12 and 16, the numbers behave like the classic 2-color van der Waerden numbers. That is, for widths $w=$ 4, 8, 12, and 16, the lengths of the smallest width-$w$ rectangle that forces an $l$-term AP of \TTs\ is four times the length of the smallest initial segment of the natural numbers that forces a monochromatic $l$-term AP under any 2-coloring. We first concentrate on the case of width 4.

\begin{figure}[h]
\begin{center}
\includegraphics[scale=0.6]{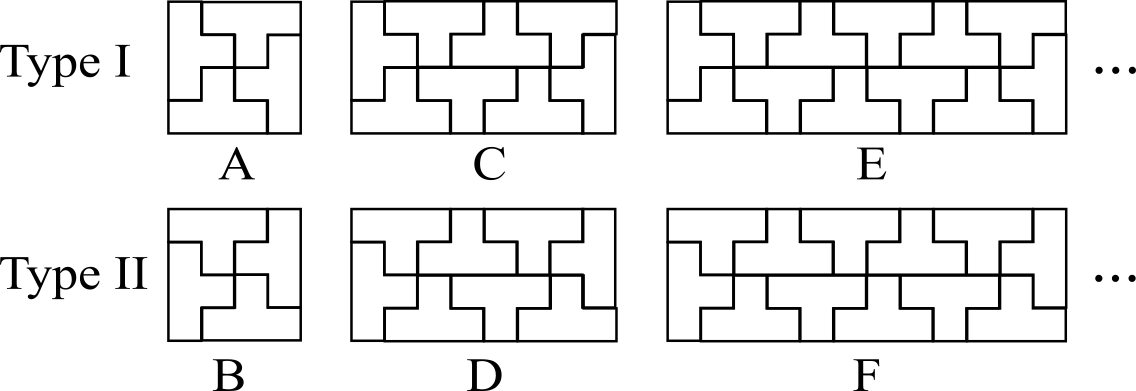}
\end{center}
\caption{Tilings of $4\times l$ rectangles are made of units like these.}
\label{fig-FellerUnits}
\end{figure}

Note that all tilings of a $4\times l$ rectangle consist of some concatenation of the units shown in Figure \ref{fig-FellerUnits}. (This follows directly from the discussion of chain graphs in Section \ref{sec-chain}.) Since long units contain long APs within themselves, we can make some easy observations: If a tiling of a $4\times l$ rectangle contains no AP of length 2, then it can contain only A units and/or B units. If such a tiling contains no AP of length 3, then it must consist of only A, B, C and/or D units. And so on. We use a two-character code to describe the orientation and vertical placement of a \TT, consisting of a letter for orientation (u, d, l, or r) and a number giving the highest row it covers, where the rows are numbered 1, 2, 3, 4 from the top down. We also number the columns $1, 2, 3, \ldots$ from left to right.

\begin{figure}[h]
\begin{center}
\includegraphics[scale=0.75]{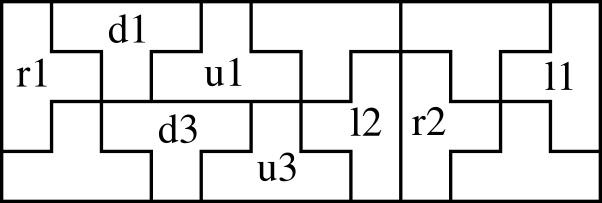}
\end{center}
\caption{Names for the orientations and vertical locations of all T-tetrominos in a $4\times N$ tiling.}
\label{fig-tetTypes}
\end{figure}

We first formalize a simple observation about arithmetic progressions.

\begin{lemma}[mod 4 Lemma]\label{lemma-mod4}
Let $A$ be an arithmetic progression of integers with at least 3 terms, and let $s$ be some integer. If, for each $a\in A$, either $a\equiv s\pmod 4$ or $a\equiv s+1\pmod 4$, then the terms of $A$ are all congruent to each other $\pmod 4$.
\end{lemma}

\noindent{\bf Proof:} It is sufficient to prove this in the case that $s=0$, since the AP can be translated by subtracting $s$ from each term without affecting the hypothesis or conclusion of the theorem. 

Suppose that $a, b$ and $c$ are three consecutive terms of the AP. Then $a+c = 2b$. If $b\equiv 0\pmod 4$, then $a$ and $c$ must both be $0\pmod 4$ in order to satisfy $a+c = 2b$. And if $b\equiv 1\pmod 4$, then $a$ and $c$ must both be $1\pmod 4$ for the same reason. In each case, the three terms must all be in the same congruence class. Since this holds for all consecutive triples, all elements of the AP must be in the same congruence class.\hfill$\Box$

Next we show that in rectangles of width 4 and APs of length $\geq 3$ it suffices to consider only type d1 \TTs.

\begin{lemma}[d1 Lemma]\label{lemma-d1}
Any tiling $\mathcal{T}$ of a width-4 rectangle contains an $l$-term AP of \TTs, for $l\geq 3$, if and only if it contains an $l$-term AP of d1 \TTs.
\end{lemma}

\noindent{\bf Proof:}  As shown in Figure \ref{fig-FellerUnits}, the only \TTs\ that occur in a $4\times N$ tiling are d1, d3, r1, r2, l1, l2, u1, and u3, and we observe (without proof) that any AP of length $\geq 3$ must consist of one of these types only. We will show that for each of those types, if some tiling contains an AP of length at least 3 of that type, then there is also an AP of d1 \TTs\ of that length in that tiling. We make recourse to Figure \ref{fig-FellerUnits} to justify the observations that follow. For the u1, r1, r2, l1 and l2 \TTs\ there will always be a d1 \TT\ joined to it in a consistent way, so that an AP in one of those types will give an AP in the joined d1 \TTs. For example, every u1 has a d1 to its immediate left in the same two rows that it occupies.

Now suppose there is an AP of u3 \TTs. Note that the u3 \TTs\ in Type-I units have their left-most square in a column number that is congruent to $1 \pmod 4$, and u3 \TTs\ in Type-II units have their left-most squares in a column number that is congruent to $2 \pmod 4$. So by the mod 4 Lemma these u3 \TTs\ will have to be all in Type-I or all in Type-II units. If it is Type-I, then there will be d1 \TTs\ above and one square to the right of each u3, giving an AP of d1 \TTs\ of the same length. And if it is Type-II, then they will be above and one square to the left. In either case we get an AP of the same length as for the u3 \TTs. Finally, we observe that a similar proof is obtained in the case of an AP of d3 \TTs.\hfill$\Box$

Every tiling $\mathcal{T}$ of a $4\times k$ rectangle can be naturally mapped to a tiling of that rectangle by A units and/or B units only. This is done by considering each unit $U$ in the tiling. $U$'s first column matches the first column of either an A unit or a B unit. If it matches the A unit, then we replace $U$ by enough A units to match the length of $U$, and if it matches the B unit, then we replace $U$ by enough B units to match the length of $U$. For example, C would be replaced by AA, and F would be replaced by BBB. We denote the resulting tiling by $AB(\mathcal{T})$.

\begin{lemma}[AB Lemma]\label{lemma-AB}
A tiling $\mathcal{T}$ of a $4\times N$ rectangle contains an $l$-term arithmetic progression of \TTs, if and only if $AB(\mathcal{T})$ contains an $l$-term arithmetic progression of \TTs.
\end{lemma}

\noindent{\bf Proof:} We note that the mapping fixes the d1 \TTs, and the result follows from the d1 Lemma.\hfill$\Box$\\

\begin{lemma}\label{lemma-4}
For all $l\geq 3$, $T_4(l) = 4W(2, l)$. 
\end{lemma}

\noindent{\bf Proof:}
First we show that $T_4(l) \geq 4W(2, l)$. Consider any two-coloring of the natural numbers of length $W(2, l)-1$ which does not contain an $l$-term monochromatic arithmetic progression.  Such a coloring exists by the definition of $W(2, l)$.  Map one of the colors to the A unit and the other color to the B unit to produce a tiling of a rectangle of size $4\times 4(W(2,l)-1)$. This tiling contains no $l$-term arithmetic progression of A units and no l-term arithmetic progression of B units.  By Lemma \ref{lemma-mod4}, this tiling contains no $l$-term arithmetic progression of d1 \TTs, and so by the d1 Lemma, contains no arithmetic progression of any orientation of \TTs.  Thus $T(4, l) \geq 4W(2,l)$.\\
%Would it be more straighforward to finish the proof with: "By Lemma 7 (AB Lemma), this tiling contains no l-term progression of T-tetrominos. Thus..." 
% Rob replies: I am not sure what this comment is suggesting. The proof is longer than I like, but not terrible, and is unambiguous... Can you flesh out the suggestion?

Next we show that $T_4(l) \leq 4W(2,l)$.
Consider a tiling $\mathcal{T}$ of a rectangle of width 4 and length $T_4(l)-4$ which has no $l$-term arithmetic progression of \TTs.  Such a rectangle exists by the definition of $T_4(l)$.  By the AB Lemma, $AB(\mathcal{T})$ contains no $l$-term arithmetic progression of \TTs, and in particular, no $l$-term arithmetic progression of A units or B units.  Map the rectangle tiled with A units and B units to a 2-coloring of the natural numbers $\{1, \ldots, T_4(l)/4 - 1\}$ such that A gets mapped to one color and B to the other. The 2-coloring resulting from the mapping contains no $l$-term monochromatic arithmetic progression, so that $W(2, l) > T_4(l)/4 - 1$, or 
$T_4(l) < 4W(2, l) + 4$, which implies $T_4(l) \leq 4W(2,l)$.\hfill$\Box$

\section{Chain Graphs and Widths 8, 12, and 16}\label{sec-chain}

\begin{figure}[h]
\begin{center}
\includegraphics[scale=0.75]{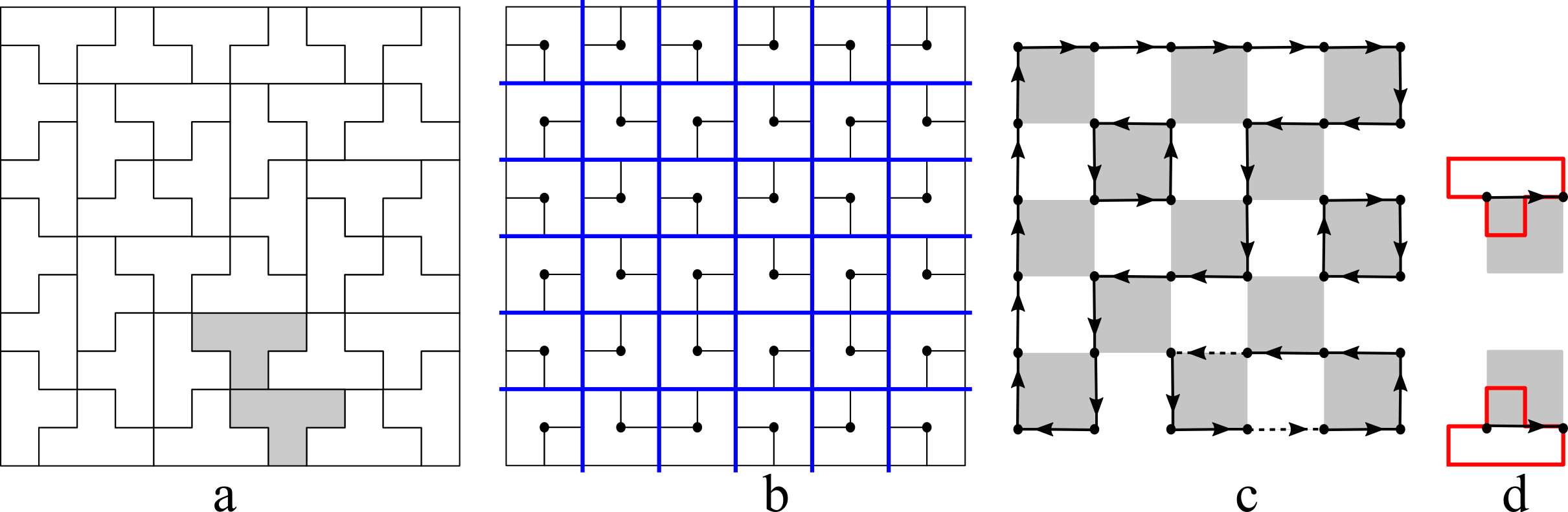}
\end{center}
\caption{A tiling (a), its $2\times 2$ blocks (b), its chain graph (c), and the correspondence between tiles and shaded arrows (d). The two shaded tiles are discussed in the proof of Lemma \ref{lemma-bijection}.}
\label{fig-chainGraph}
\end{figure}

To address widths 8, 12 and 16 we introduce some notions from a paper of Korn and Pak \cite{KornPak}. Every tiling (Figure \ref{fig-chainGraph}a) gives rise to a directed graph as follows: Subdivide the rectangle into $2\times 2$ blocks, as shown in Figure \ref{fig-chainGraph}b. Korn and Pak show, based on Walkup's paper, that each \TT\ must have three of its squares in one block (its ``majority'' block) and one square in an adjacent block  (its ``minority'' block). Place a vertex into the center of each block, and for each \TT, draw a directed edge from the vertex in its majority block to the vertex in its minority block, as shown in Figure \ref{fig-chainGraph}c. They call the resulting digraph the {\em chain graph} for that tiling. Korn and Pak call the $2\times 2$ boxes inside the chain graph {\em antiblocks}, on which they impose a gray/white checkerboard coloring with gray in the corners, as shown in Figure \ref{fig-chainGraph}c. Each arrow has a gray block on its left or its right, as seen when walking along that arrow in its direction. See Figure \ref{fig-chainGraph}d. We define an AP of {\em shaded arrows} to be a sequence of equally-spaced arrows, pointing in the same direction, having the same left/right shading.

\begin{lemma}\label{lemma-bijection}
For any $l\geq 3$, a tiling has an AP of \TTs\ of length $l$ if and only if its chain graph has an AP of shaded arrows of length $l$.
\end{lemma}

\noindent{\bf Proof:}
It follows from the discussion in the proof of Lemma 5 in \cite{KornPak} that every right-pointing shaded arrow corresponds to a tile exactly as shown in Figure \ref{fig-chainGraph}d, depending on whether its left or right side is shaded. Those figures can be rotated to show the correspondence for the other three orientations. Thus any $l$-term AP of shaded arrows gives rise to a $l$-term AP in the associated \TTs. 

The other direction is slightly more subtle. For example, the two shaded \TTs\ in Figure \ref{fig-chainGraph}a form a 2-term AP, but their corresponding shaded arrows (dashed in the Figure \ref{fig-chainGraph}c) do not. To show that this cannot happen with APs of length 3 or greater, we need another lemma.

\begin{lemma}[The dx/dy Lemma]\label{lemma-dxdy}
Suppose we are given an AP of \TTs\ of length 3 or greater, and let the pair $(dx, dy)$ give the translation between consecutive tiles in the AP. Then $(dx, dy)$ is congruent to either $(0, 0)$ or $(2, 2)\pmod 4$.
\end{lemma}

\begin{figure}[h]
\begin{center}
\includegraphics[scale=0.4]{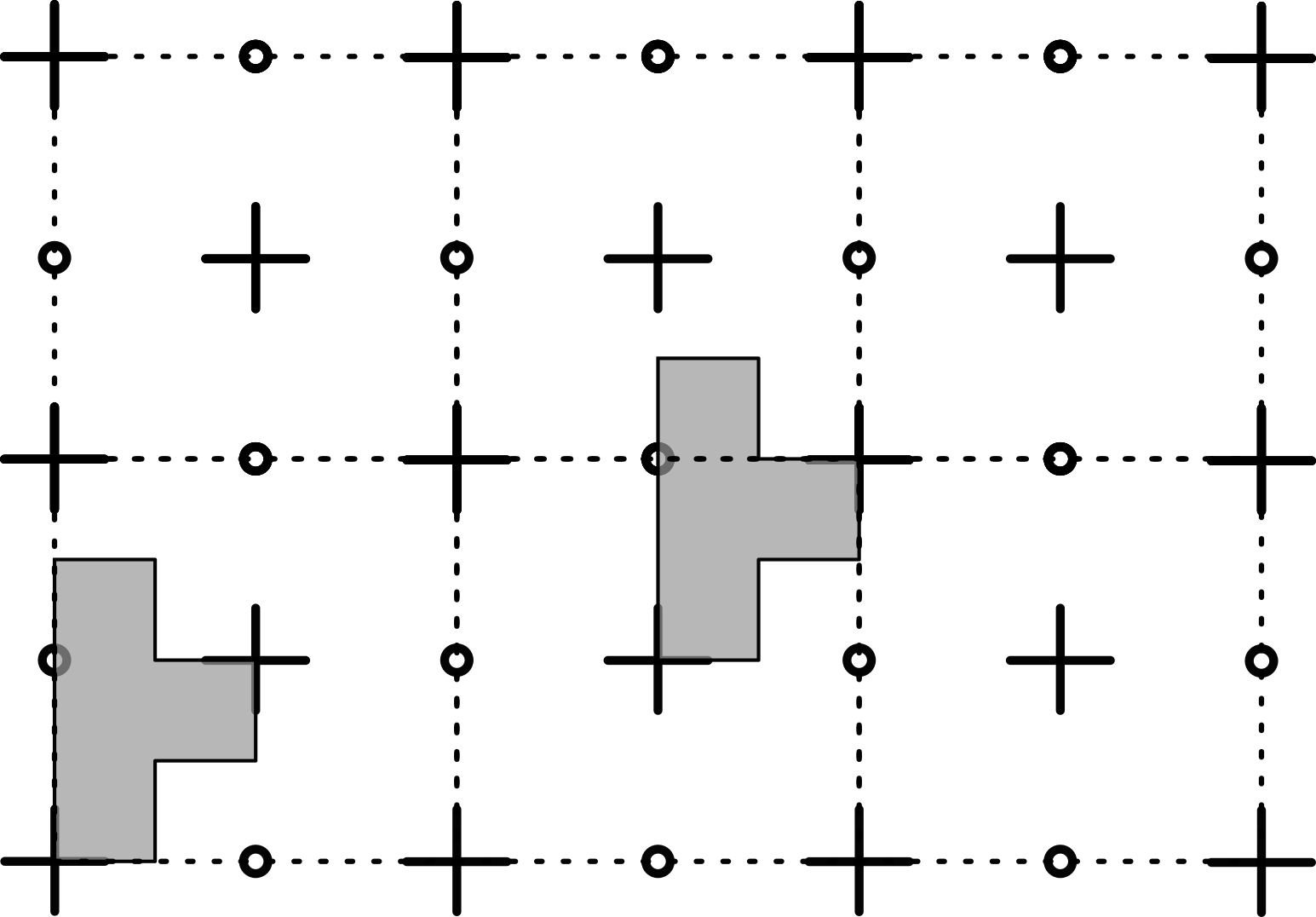}
\end{center}
\caption{The cut and cornerless points.}
\label{fig-cutCornerless}
\end{figure}

\noindent{\bf Proof:}
Walkup's proof makes use of {\em cut} and {\em cornerless} points, which are, respectively, grid points where four tile corners must meet, and grid points where no tile corner may lie, in any tiling of a rectangle. Figure \ref{fig-cutCornerless} shows cut points as crosses and cornerless points as circles, in $4\times 4$ blocks that fill the rectangle periodically. More precisely, if the lower-left corner of the rectangle has coordinates $(0, 0)$, then the cut and cornerless points are exactly those points $(x, y)$ where both $x$ and $y$ are even; and $x+y\equiv 0\pmod 4$ for cut points, and $x+y\equiv 2\pmod 4$ for cornerless points. 
%MMM How is this ^^^ for being precise? Yes!
The two shaded tiles in that figure show the only two possible placements of a \TT, up to rotation or reflection with respect to those $4\times 4$ blocks. Consider the \TT\ on the bottom left. Of the sixteen possible pairs $(dx, dy) \pmod 4$, twelve can be immediately ruled out by considering cut and cornerless points. For example, $(dx, dy) = (1, 0)$ is ruled out, because that tile, slid one square to the right, results in a tile with a corner at a cornerless point, as well as cutting two arms of a cross. Only $(0, 0), (0, 1), (2, 2)$ and $(2, 3)$ permit a valid translation. But as $(0, 1)$ and $(2,3)$ do not permit a {\em second} translation, they can be ruled out in any AP of length 3 or greater. The proof is similar for the other starting position.\hfill$\Box$

We may now finish the proof of Lemma \ref{lemma-bijection}. Suppose we have an AP of at least three \TTs\, and consider the associated shaded arrows. By the dx/dy Lemma, their corresponding arrows must have their shading on the same side, and therefore have the same orientation, giving an AP of shaded arrows of the same length.\hfill$\Box$
%Because we rearrenged the paper, we have not yet started the proof of Theroem 12.  We could remove the first sentence of the paregraph, the begin the second sentence with "Lemma 11 follows from..."
% RRR Rob replies: This seems to be based on a previous version...

\begin{figure}[h]
\begin{center}
\includegraphics[scale=0.7]{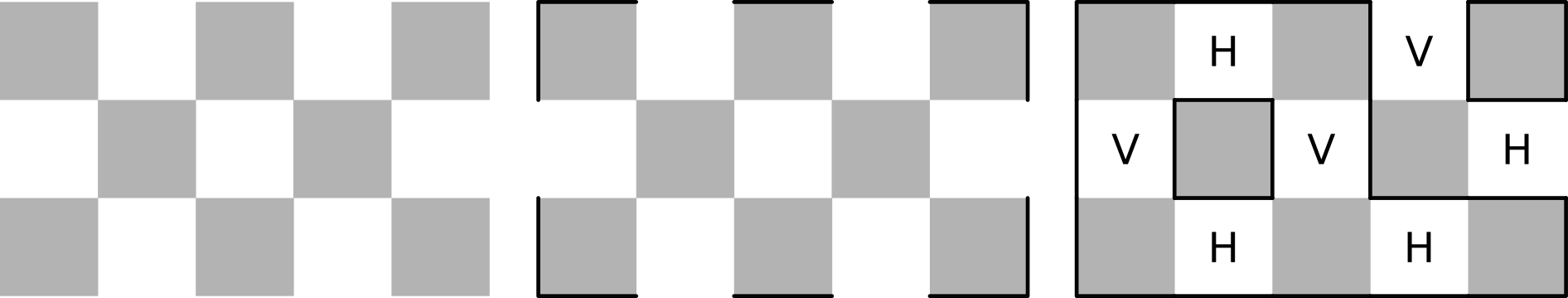}
\end{center}
\caption{Illustration of steps 1 and 2 of the HV Lemma.}
\label{fig-hvtheorem}
\end{figure}

We give one more lemma to help us finish the proof of Theorem \ref{thm-4812} for widths 8, 12 and 16. This follows from the discussion in Section 4 of \cite{KornPak}, so we do not prove it here. It is illustrated in Figure \ref{fig-hvtheorem}.

\begin{lemma}[The HV Lemma]\label{lemma-hv}
Every chain graph can be constructed as follows from the gray/white antiblocks: 1. Add a chain edge along every side of a gray antiblock that does not have an adjacent white antiblock; 2. For each white antiblock, add either both horizontal edges, or both vertical edges; 3. Arbitrarily orient each resulting cycle.\hfill$\Box$
\end{lemma}

\begin{thm}\label{thm-4812}
For all $l\geq 3$, $T_4(l) = T_8(l) = T_{12}(l) = T_{16}(l) = 4W(2, l)$. 
\end{thm}

\noindent{\bf Proof:} The case for $T_4(l)$ was proved in Lemma \ref{lemma-4}. This section deals with widths 8, 12, and 16. 

Observe that for any width $w\geq 4, T_w(l) \leq 4W(2, l)$. For, given a tiling of a $w\times 4W(2, l)$ rectangle, the HV Lemma says that each gray anti-block in the top row yields an arrow on its upper edge, and we may consider the left/right orientation of those arrows to be a 2-coloring of the set $\{1, 2, \ldots, W(2, l)\}$. By the definition of $W(2, l)$ this will contain an AP of length $l$ of shaded arrows, and hence of \TTs. (Note that this gives the desired improvement of the bound given in the proof of Theorem \ref{thm-boundary}.)

The proof of Lemma \ref{lemma-4} showed that $T_4(l) > 4W(2, l)-4$, by constructing a row of A blocks and B blocks. To show that $T_8(l) > 4W(2, l)-4$ we can take two copies of such a row, one placed above the other, each of which contains no AP of length $l$. Suppose this tiling contained an AP of \TTs\ of length $l$ and let $(dx, dy)$ be its translation. We cannot have $dy=0$, for then the AP would lie in a single row whose A/B blocks were constructed to contain no such AP. The structure of A/B blocks does not permit $dy = 2$. And if $dy\geq 4$ there is not enough vertical space to fit the AP, since $l\geq 3$. This completes the proof for width $w=8$.

\begin{figure}[h]
\begin{center}
\includegraphics[scale=1.1]{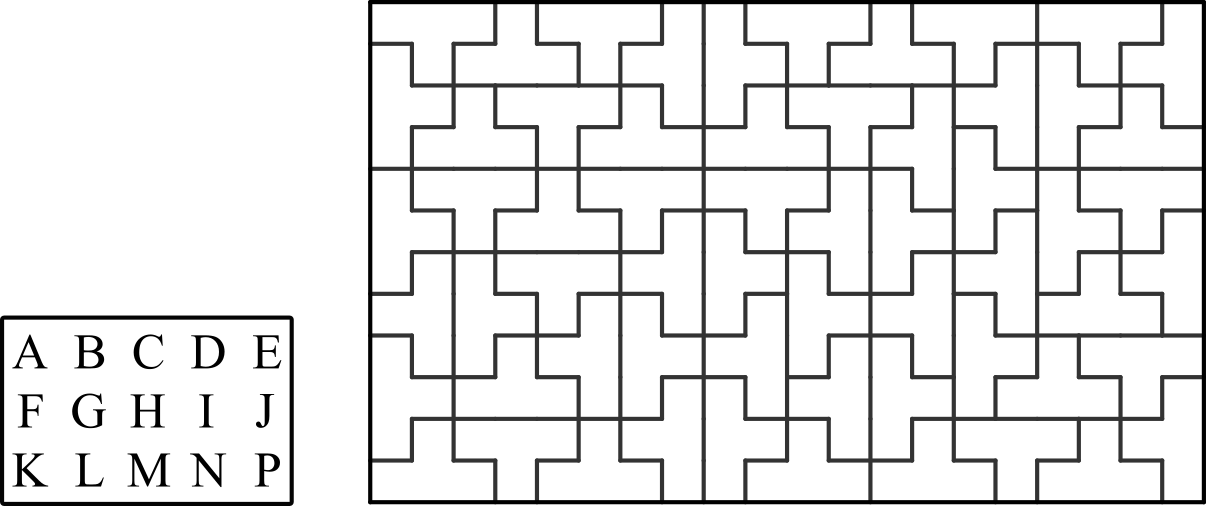}
\end{center}
\caption{Every 2-coloring of a $3\times 5$ grid has a 3-term monochromatic AP, but not every tiling of a $12\times 20$ grid.}
\label{fig-apon12x20}
\end{figure}

For widths 12 and 16 the situation is different. We can not use A/B blocks, for example, to show that $T_{12}(3) \geq 4W(2, 3) = 36.$ This is because any 2-coloring of the squares of a $3\times 5$ grid must contain a monochromatic 3-term AP, as we will now show by considering the $3\times 5$ rectangle on the left in Figure \ref{fig-apon12x20}. Suppose there was a red/blue-coloring of those 15 squares with no 3-term AP. We may assume that square H is colored red. Then at least one of A or P must be blue, and at least one of E or K must also be blue. This gives two cases up to symmetry. In case 1, we assume that A and E are blue. Then C must be red, which forces M to be blue, which forces both G and I to be red, giving a contradiction. The other case is similar. This implies that any tiling of a $12\times 20$ rectangle with A/B blocks must contain a 3-term AP, yet there does exist a tiling of a $12\times 20$ rectangle, as shown on the right in that figure, containing no 3-term AP.

\begin{figure}[h]
\begin{center}
\includegraphics[scale=0.32]{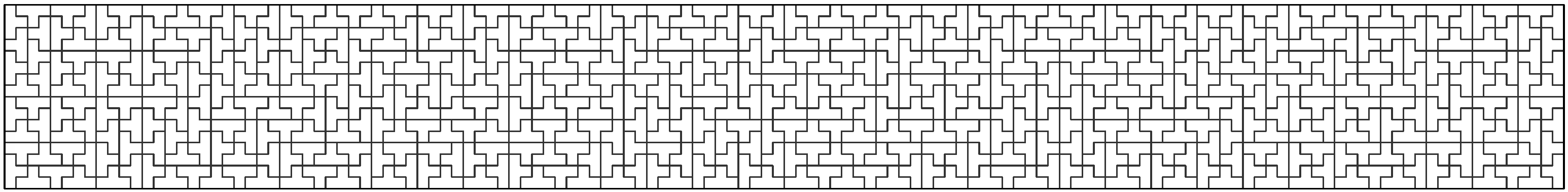}
\end{center}
\caption{A tiling of a $16\times 136$ rectangle with no AP of length 4.}
\label{fig-16x136}
\end{figure}

So instead of using just A/B blocks, for width $w=12$ we have constructed tilings of all $12\times 4b$ rectangles containing no AP of length 3, for $4\leq 4b\leq 32 = 4W(2, 3) - 4$. (All tilings that we have constructed can be found in the supplementary materials \cite{supp}.)  This proves that $T_{12}(l)\geq 4W(2, l)$ when $l=3$. For $l>3$, we may use stacked A/B tilings as we did for width 8, but with three identical rows instead of two. Finally, for width $w=16$, we have constructed tilings of all $16\times 4b$ rectangles containing no AP of length 3, for $4\leq 4b\leq 32 = 4W(2, 3) - 4$, as well as tilings of all $16\times 4b$ rectangles containing no AP of length 4, for $4\leq 4b\leq 136 = 4W(2, 4) - 4$. (See Figure \ref{fig-16x136}.) For $l>4$ we may use four stacked identical A/B tilings.\hfill$\Box$

By contrast we note that while $T_{20}(4) = 140 = 4W(2, 4)$, $T_{20}(3) = 24$, which is less than $36 = 4W(2, 3)$. See Figure \ref{fig-smallNumbers}. Theorem \ref{thm-4812} may be generalized to arbitrary widths, as long as the lower bound on $l$ is made suitably large, since for any width, a large enough $l$ will permit only horizontal APs, and stacked A/B tilings will give optimal tilings.
% MMM Is this ^^^ a good "by contrast?" as foreshadowed at the end of the Introduction? Yes, I think that connects well.  I did change the 136 to 140. 
% RRR Good catch, Emily!

\section {Exact values}
\begin{figure}[h]
\begin{center}
\includegraphics[scale=0.45]{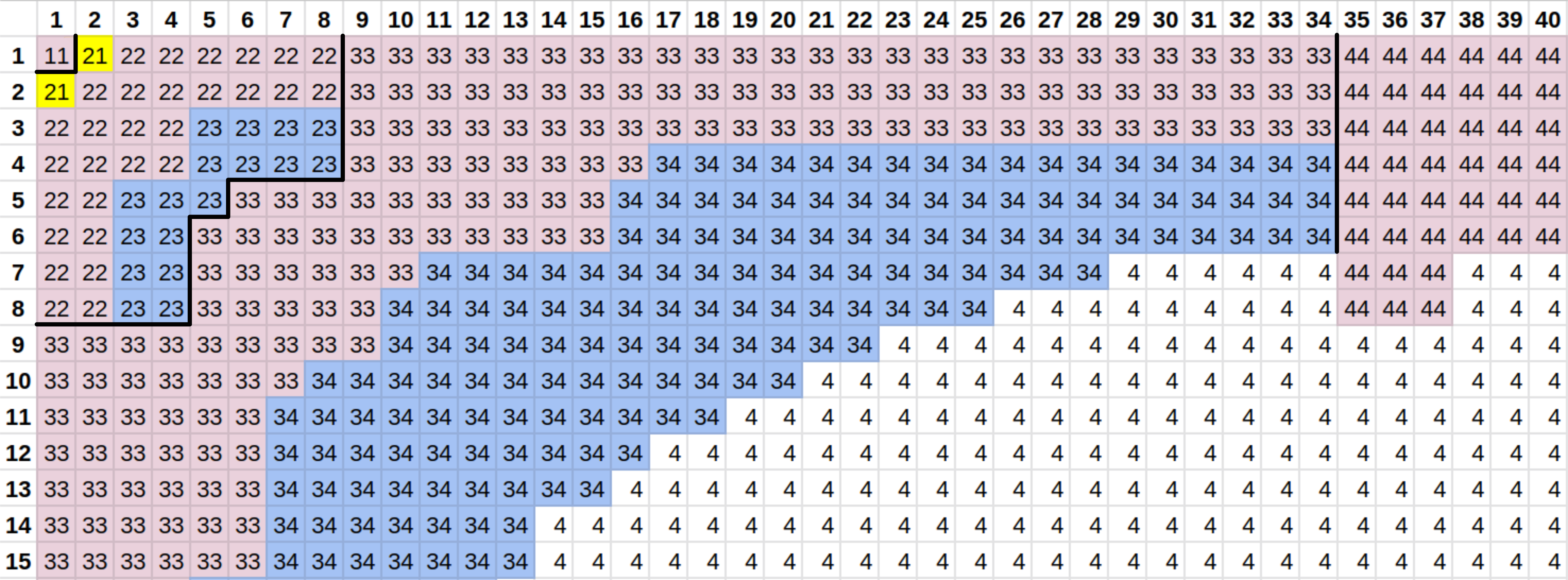}
\end{center}
\caption{Some exact values of $L(h, w)$ and $L_{vdW}(h, w)$. The first digit in row $h$ column $w$ is $L(4h, 4w)$, and the second digit is $L_{vdW}(h, w)$. Color indicates the difference between the digits. Where no first digit is given, we do not know the corresponding value of $L$. The bold lines separate different values of $L$.}
\label{fig-smallNumbers}
\end{figure}

Let us define $L(h, w)$ to be the greatest value $l$ such that $(h, w)\rightarrow l$. Figure \ref{fig-smallNumbers} shows all values of $L(h, w)$ that we have been able to compute so far, plus those given in the proof of Theorem \ref{thm-4812}. 
One consequence of the AB Lemma, and Theorem \ref{thm-4812}, is that for $h$ in $\{4, 8, 12, 16\}$ if $w_1 \geq w_2$, then $L(h, w_1) \geq L(h, w_2)$. So far, this is the only monotonicity result that we have been able to prove, as well as that suggested by the generalization mentioned in the previous paragraph. And, without the assumption of monotonicity, every single proof that $L(h, w) = l$ requires the construction of a new tiling of an $h\times w$ rectangle with no AP of length $l+1$, even if we have such a tiling of a larger rectangle; and also a proof that every tiling of an $h\times w$ rectangle must contain an AP of length $l$, even if we have proved that for some smaller rectangle. Most of our results have been found using SAT solvers, that is, computer programs for solving problems represented as boolean formulae.

\begin{figure}[h]
\begin{center}
\includegraphics[scale=0.6]{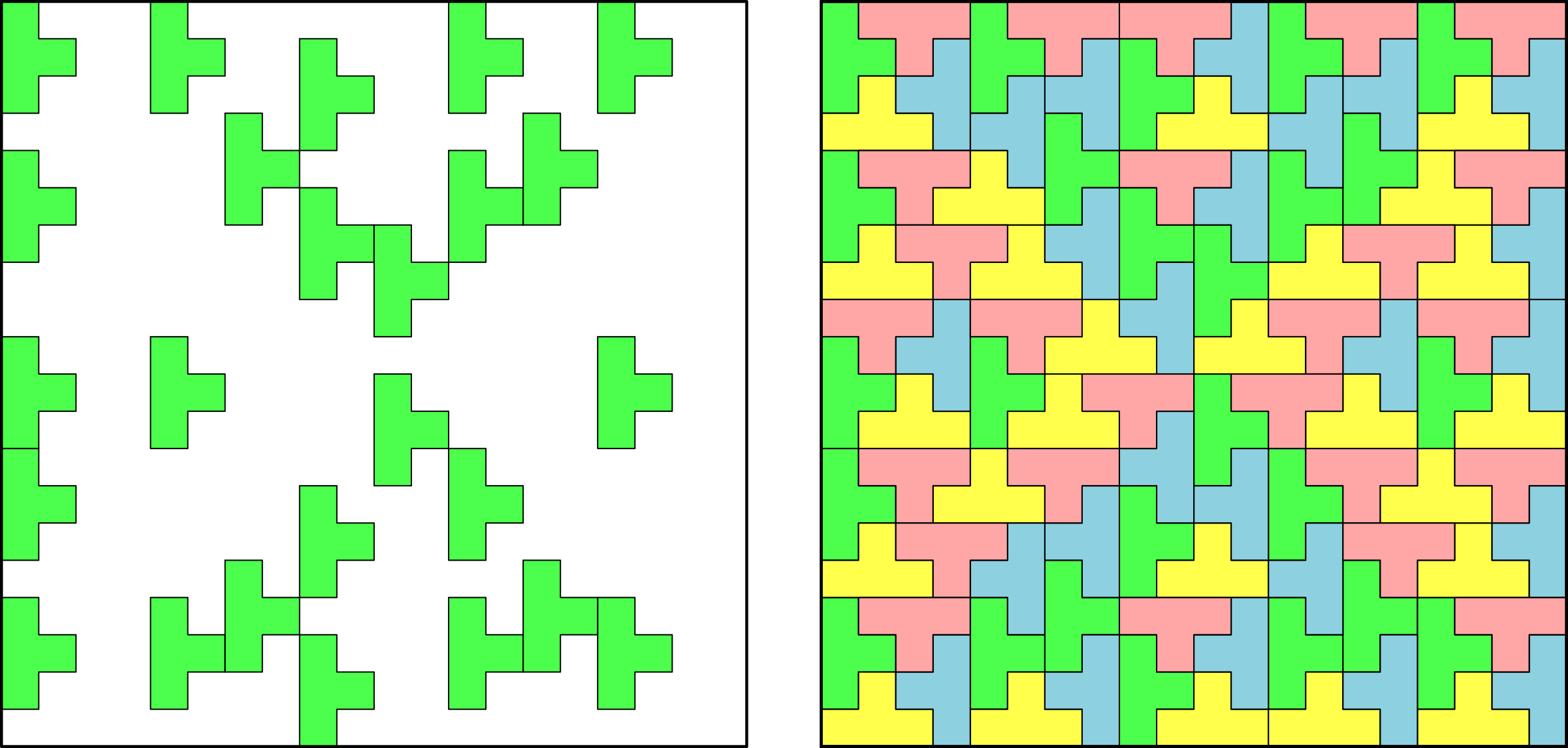}
\end{center}
\caption{A rotationally symmetric $20\times20$ tiling with no 3-term AP.}
\label{fig-rotSymmetry}
\end{figure}

Code for generating SAT instances (boolean formulas representing a tiling problem) and solving them using pysat can be found in the supplementary materials \cite{supp}. Note that some values took several days to compute. Sometimes searching for symmetric tilings yielded a speedup, as well as some beautiful tilings. For example, Figure \ref{fig-rotSymmetry} gives a rotationally-symmetric maximal rectangle with no 3-term AP.

\subsection{2D van der Waerden Numbers}
There seems to be no literature on APs of colored squares in a colored rectangular grid, so we mention it briefly here. As mentioned in the discussion of Figure \ref{fig-apon12x20}, every 2-coloring of the squares of a $3\times 5$ grid must contain a monochromatic 3-term AP. A coloring with alternating vertical stripes shows that a 4-term AP may be avoided. If we define $l = L_{vdW}(h, w)$ to be the greatest value $l$ such that every 2-coloring of an $h\times w$ grid must contain a monochromatic AP of length $l$, then we've shown that $L_{vdW}(3, 5) = 3$. Because of the connection between $L$ and $L_{vdW}$ given by the AB Lemma we wished to compute more values of the latter. Thanks to the obvious monotonicity of $L_{vdW}$ this was much easier, and the results are shown in Figure \ref{fig-smallNumbers}. That figure shows all values where $L_{vdW}(h, w) = 3$, and the following list of $(h, w)$ pairs show all sizes we know where $L_{vdW}(h, w) = 4$: $(4, 100), (10, 50), (15, 40), (20, 32), (24, 24)$, together with those implied by symmetry ($L_{vdW}(h,w) = L_{vdW}(w,h)$) and monotonicity. To prove that $L_{vdW}(h, w) = 4$, we used SAT solvers to show that no 2-coloring of an $h\times w$ rectangle without a monochromatic 4-term AP existed, and then to find a 2-coloring of an $h\times w$ rectangle that did not contain a monochromatic 5-term AP.

\subsection{Further Investigations} We conclude with two conjectures and a suggested research direction:
\begin{itemize}
    \item {\bf Conjecture:} The function $L(h, w)$ is monotonic. 
    \item {\bf Conjecture:} The inequality $\max |L(4h, 4w)-L_{vdW}(h, w)| \leq 1$ holds for all $h$ and $w$, as suggested by the limited data shown in Figure \ref{fig-smallNumbers}.
    \item Find tight upper and lower bounds on $L(h, h)$, the diagonal entries in Figure \ref{fig-smallNumbers}.
\end{itemize}

\section{Acknowledgements}
The authors thank Eli Paul and other students at the Hampshire College Summer Studies in Mathematics for help computing some values of $L_{vdW}(h, w)$, and the NSF-sponsored CREU and CRA-W for supporting the first author.

\bibliographystyle{elsarticle-num} 
\bibliography{TileRamseyNumbers}

\begin{thebibliography}{1}
\expandafter\ifx\csname url\endcsname\relax
  \def\url#1{\texttt{#1}}\fi
\expandafter\ifx\csname urlprefix\endcsname\relax\def\urlprefix{URL }\fi
\expandafter\ifx\csname href\endcsname\relax
  \def\href#1#2{#2} \def\path#1{#1}\fi

\bibitem{Walkup}
D.~W. Walkup, Covering a rectangle with {$T$}-tetrominoes, Amer. Math. Monthly
  72 (1965) 986--988.

\bibitem{ReidKlarner}
M.~Reid, \href{http://dx.doi.org/10.1016/j.jcta.2004.10.010}{Klarner systems
  and tiling boxes with polyominoes}, J. Combin. Theory Ser. A 111~(1) (2005)
  89--105.
\newblock \href {http://dx.doi.org/10.1016/j.jcta.2004.10.010}
  {\path{doi:10.1016/j.jcta.2004.10.010}}.
\newline\urlprefix\url{http://dx.doi.org/10.1016/j.jcta.2004.10.010}

\bibitem{KornPak}
M.~Korn, I.~Pak, \href{http://dx.doi.org/10.1016/j.tcs.2004.02.023}{Tilings of
  rectangles with {T}-tetrominoes}, Theoret. Comput. Sci. 319~(1-3) (2004)
  3--27.
\newblock \href {http://dx.doi.org/10.1016/j.tcs.2004.02.023}
  {\path{doi:10.1016/j.tcs.2004.02.023}}.
\newline\urlprefix\url{http://dx.doi.org/10.1016/j.tcs.2004.02.023}

\bibitem{Merino}
C.~Merino, On the number of tilings of the rectangular board with
  {T}-tetrominoes, Australas. J. Combin. 41 (2008) 107--114.

\bibitem{Goddard}
W.~Goddard, Almost tilings with t-tetrominos, presented at the Cumberland
  Conference, write-up unavailable (May 2007).

\bibitem{Hochberg}
R.~Hochberg,
  \href{http://www.sciencedirect.com/science/article/pii/S0012365X14003525}{The
  gap number of the t-tetromino}, Discrete Mathematics 338~(1) (2015) 130 --
  138.
\newblock \href {http://dx.doi.org/https://doi.org/10.1016/j.disc.2014.09.001}
  {\path{doi:https://doi.org/10.1016/j.disc.2014.09.001}}.
\newline\urlprefix\url{http://www.sciencedirect.com/science/article/pii/S0012365X14003525}

\bibitem{Zhan}
S.~Zhan, Tiling a deficient rectangle with t-tetrominoes, preprint (August
  2012).

\bibitem{GRS}
R.~L. Graham, B.~L. Rothschild, J.~H. Spencer, Ramsey theory, Wiley Series in
  Discrete Mathematics and Optimization, John Wiley \& Sons, Inc., Hoboken, NJ,
  2013, paperback edition of the second (1990) edition [MR1044995].

\bibitem{supp}
Feller, Hochberg,
  \href{https://github.com/RobertHochberg/FellerHochberg}{Supplementary
  materials} (2024).
\newline\urlprefix\url{https://github.com/RobertHochberg/FellerHochberg}

\end{thebibliography}

\end{document}